\documentclass[fleqn]{mat01}
\usepackage{times,mathtimy,amssymb,latexsym}
\begin{document}

\setcounter{page}{317} \firstpage{317}


\def\pa{\S\kern.15em}

\def\theor{\trivlist \item[\hskip \labelsep{\bf Theorem.}]}
\def\prob{\trivlist \item[\hskip \labelsep{\it Problem.}]}

\def\({{\rm (}}
\def\){{\rm )}}
\def\[{{\rm [}}
\def\]{{\rm ]}}
\def\ba{\backslash}
\def\1{\hbox{\small 1}\!\!\hbox{\normalsize 1}}
\def\n{\not=\!\!\!\not=}
\def\SL{\textrm{SL}}
\def\Gal{\textrm{Gal}}
\def\ind{\textrm{ind}}
\def\GL{\textrm{GL}}

\title{Equivariant embeddings of Hermitian symmetric spaces}

\markboth{L~Clozel}{Equivariant embeddings of Hermitian symmetric spaces}

\author{L~CLOZEL}

\address{Universit\'{e} Paris-Sud Orsay, Math\'{e}matiques B$\hat{\rm a}$t 425, 91405
Orsay, France\\
\noindent E-mail: Laurent.Clozel@math.u-psud.fr}

\volume{117}

\mon{August}

\parts{3}

\pubyear{2007}

\Date{MS received 28 November 2006; revised 3 April 2007}

\begin{abstract}
We prove that equivariant, holomorphic embeddings of Hermitian
symmetric spaces are totally geodesic (when the image is not of
exceptional type).
\end{abstract}


\maketitle

\section{Introduction}

Let $H,G$ be connected semi-simple Lie groups and $X_{H}$, $X_{G}$ the associated
symmetric spaces. We assume that they are Hermitian. An equivariant embedding is a
pair $(F,f)$ where $F\!\!:H\rightarrow G$ is a homomorphism, $f\!\!:X_{H}\rightarrow
X_{G}$ is a holomorphic map and
\begin{equation*}
f(h\cdot x)=F(h)f(x), \quad x\in X_{H},h\in H.
\end{equation*}

We assume that $H,G$ have no compact factors and that $f$ is injective. Then, as is
easily checked, the kernel of $F$ is finite. Replacing $H$ by its image, we will also
assume $F$ injective and therefore identify $H$ with its image in $G$.

Such maps have been classified by Satake \cite{8} and Ihara \cite{3} when $X_{H}$ is
{\it totally geodesic} in $X_{G}$. The purpose of this note is to show the following
theorem.

\begin{theor}{\it
Assume $G$ has no factors of exceptional type. Then any
equivariant embedding $X_{H}\rightarrow X_{G}$ is totally
geodesic.}\vspace{.5pc}
\end{theor}

We should emphasize the rather surprising content of this result when compared with
the case of compact Hermitian symmetric spaces. If $G$ is compact, the symmetric space
$X_{G}$ (assumed Hermitian) is a generalized Grassmanian. The natural maps of
algebraic geometry between Grassmanians~--~in particular the Veronese and Segre
embeddings~--~are holomorphic and equivariant with respect to natural maps of the
associated groups. Very few are totally geodesic: in fact by duality between compact
and non-compact symmetric spaces, the {\it totally geodesic} equivariant maps between
compact spaces correspond to those between non-compact spaces, which are quite rare
(see~\cite{2}). However, this result becomes more natural from the `global' point of
view, i.e., if one considers arithmetic quotients of the symmetric spaces.

Assume $H,G$ are semi-simple groups defined over $\mathbb{Q}$, $F\!\!:H\rightarrow G$
is defined over $\mathbb{Q}$ and $f\!\!:X_{H}\rightarrow X_{G}$ is an equivariant
embedding. For suitable arithmetic subgroups $\Delta\subset H(\mathbb{Q})$ and
$\Gamma\subset G(\mathbb{Q})$, $f$ defines a holomorphic map
\begin{equation*}
g\!\!:S_{H}\rightarrow S_{G},
\end{equation*}
where $S_{H}=\Delta\ba X_{H}$ and $S_{G}=\Gamma\ba
X_{G}$.\pagebreak

In this situation, recall that $S_{H}$ and $S_{G}$ have a
remarkable family of distinguished points, the CM-points or
special points \cite{5}. Also note that $S_{H}$, $S_{G}$ are in
fact algebraic varieties over $\mathbb{C}$, and that $g(S_{H})$ is
an algebraic subvariety of $S_{G}$ by a theorem of Borel. Assume
$g(S_{H})$ has {\it one} CM-point. By using the action of
$H(\mathbb{Q})$ on $X_{H}$ one easily sees that it has a {\it
dense} subset of CM-points for the complex topology. A~conjecture
of Andr\'e \cite{1} and Oort \cite{7} then implies that $g(S_{H})$
is a  totally geodesic submanifold of $S_{G}$ (and $X_{H}$ is a
totally geodesic submanifold of $X_{G}$).

It is not obvious that $g(S_{H})$ should have one CM-point; note, however, the
following. The Hermitian symmetric spaces are open subspaces of their compact
duals~--~generalized Grassmanians. An equivariant holomorphic embedding will generally
be given by a natural holomorphic map between the compact duals. Given the
$\mathbb{Q}$-structure, CM-points correspond to subspaces (in the Grassmanians)
verifying some rationality conditions. It is natural to expect these to be preserved.
The embedding of the symmetric space for $SU(p,1)$, $X_{p,1}$, into $X_{P,Q}$ where
$P=\binom{p}{k}$, $Q=\left(\begin{smallmatrix}p
\\k-1\end{smallmatrix}\right)$ \cite{8} gives a very graphic example.

Another strong motivation for the theorem is given by Mok's rigidity results. Assume
for simplicity that $H$ is irreducible over $\mathbb{Q}$ and $\hbox{rk}(H)>1$ (this is
the {\it real} rank). Then Mok (Ch.~6, Thm~4.1 of \cite{4})~--~see also the discussion
at the beginning of ch.~9~--~has shown that {\it any} holomorphic map
$S_{H}\rightarrow S_{G}$ is totally geodesic. If $F\!\!:H(\mathbb{R})\rightarrow
G(\mathbb{R})$ (we now denote the Lie groups by $G(\mathbb{R)}$, $H(\mathbb{R})$ as we
will be using rationality arguments) is given and if $F$ is $G(\mathbb{R})$-conjugate
to a map defined over $\mathbb{Q}$, Mok's theorem implies our local assertion. More
generally, assume $F$ given, and assume that there exists a totally real number field
$L$ and a map $F_{L}\!\!:H\rightarrow G$ defined over $L$ such that, for each real
prime $v$ of $L$ (thus $L_{v}\cong\mathbb{R})$,
\begin{equation*}
F_{L,v}\!\!:H(\mathbb{R})\rightarrow G(\mathbb{R})
\end{equation*}
is conjugate to $F$. Then, again using Mok's results, we deduce that $F$ is totally
geodesic.

The set of homomorphisms $F\!\!:H\rightarrow G$, over an
algebraically closed field, and modulo $G$-conjugation, is {\it
discrete} (homomorphism of semi-simple groups up to conjugacy are
rigid). Thus $F\!\!:H(\mathbb{R})\rightarrow G(\mathbb{R})$ is
$G(\mathbb{R})$-conjugate to a map  $F_{1}$ defined  over
$\bar{\mathbb{Q}}$; the $G$-conjugacy class of $F_{1}$ is an
irreducible variety. If it is defined over $\mathbb{Q}$, a theorem
of Moret--Bailly \cite{6} implies that there is a totally real
number field $L$, and a map
$F_{L}\!\!:H\otimes_{\mathbb{Q}}L\rightarrow
G\otimes_{\mathbb{Q}}L$ verifying our condition.

It is of course, difficult to compute the field of rationality of the class associated
to $F$. One may, however, pose the following:

\begin{prob}
If $H,G$ be semisimple groups over $\mathbb{Q}$ and
$F\!\!:H\rightarrow G$ a homomorphism defined over $\mathbb{R}$,
does there exist a totally real field $L$ and
$F_{L}\!\!:H\rightarrow G/L$ such that $F_{v}$ is
$G(\mathbb{R})$-conjugate to $F$ at each real prime of
$L$?\vspace{.5pc}
\end{prob}

Finally, Mok has informed us that he could prove the theorem even for exceptional $G$.
His proof, however, is more difficult and necessitates global geometric computations.

\section{Reductions}

Let $G$ be a connected semi-simple Lie group, with finite center and no compact
factor, associated to a Hermitian symmetric space $X$. Fix a point $x\in X$. Then $x$
defines a maximal compact subgroup $K\subset G$ and a Cartan involution $\theta$ on
$\mathfrak{g}=\mathrm{Lie}(G)$. Let
\begin{equation*}
\mathfrak{g}=\mathfrak{k}\oplus\mathfrak{p}
\end{equation*}
be the Cartan decomposition. There exists an element $\zeta\in Z(K)$ such that
$\hbox{Ad}(\zeta)$ induces on $\mathfrak{p}$ the multiplication by $i=\sqrt{-1}$
defining the complex structure. Then $\zeta^{2}\in Z(K)$ induces, by the adjoint
action, the Cartan involution. By construction this holomorphic structure is
$G$-equivariant: if $x'=g\cdot x$ the associated data are obtained by conjugation by
$g$. In particular, $\zeta'=\hbox{Ad}(g)\zeta\in K'$ is well-defined by $x'$ since
$\zeta$ is $K$-invariant, and this family of quasi-complex structures defines the
holomorphic structure on $X$.

Now assume $H\subset G$, $f\!\!:X_{H}\rightarrow X_{G}$ verify our conditions. Fix a
base point $x\in X_{H}$. This defines maximal compact subgroups $K_{H}\subset K_{G}$.
(We will drop indexes for the group $G$). Thus
\begin{align*}
\mathfrak{g}&=\mathfrak{k}\oplus\mathfrak{p}\ ,\\[.2pc]
\mathfrak{h}&=\mathfrak{k}_{H}\oplus\mathfrak{p}_{H}
\end{align*}
and the (injective) map $F\!\!:\mathfrak{h}\rightarrow\mathfrak{g}$ has the following
properties:
\begin{align}
&F(\mathfrak{k}_{H})\subset \mathfrak{k},\\[.2pc]
&F(X)=F_{c}(X)+F_{p}(X),\\[.2pc]
&(X\in\mathfrak{p}_{H},F_{c}(X)\in \mathfrak{k}\ ,\
F_{p}(X)\in\mathfrak{p})\nonumber\\[.2pc]
&F_{p}(\iota_{H}X)=\iota_{G}\ F_{p}(X),
\end{align}
where $\iota_{H}$, $\iota_{G}$ are `multiplication by $\sqrt{-1}$' on
$\mathfrak{p}_{H}$, $\mathfrak{p}$, given by $\zeta_{H}$, $\zeta_{G}$. Conversely, if
a morphism $F\!\!:\mathfrak{h}\rightarrow\mathfrak{g}$ verifies (1)--(3), $F$ defines
a map $H/K_{H}\overset f\rightarrow G/K_{G}$, holomorphic at $x=eK_{H}$ and in fact at
every point by a computation similar to that as above. Note that $f$ is a totally
geodesic immersion if and only if,
\begin{equation*}
F(\mathfrak{p}_{H})\subset \mathfrak{p}\ ,\ \hbox{i.e.}, \ \mathrm{if}\ F_{c}\equiv 0.
\end{equation*}
(see p.~47 ff. of \cite{8})

Identifying $\mathfrak{k}_{H}$ with a subalgebra of $\mathfrak{k}$ by (1), we note
that the two components $F_{c}$ and $F_{p}$ are $\mathfrak{k}_{H}$-equivariant.
Moreover, let $\mathfrak{h}=\oplus\mathfrak{h}_{i}$ be a decomposition of
$\mathfrak{h}$ in simple factors. Then $\zeta_{H}$ or $\iota_{H}$ decomposes
accordingly, so the restriction $F_{i}$ to $\mathfrak{h}_{i}$ again verifies the
conditions. Thus we may assume that $\mathfrak{h}$ is simple.

In this case it is known (see S Helgason, Differential Geometry
and Symmetric Spaces, ch.~VIII, \pa\!5) that the (real)
representation of $\mathfrak{k}_{H}$ on $\mathfrak{p}_{H}$ is
irreducible. The $\mathfrak{k}_{H}$-map
$F_{c}\!\!:\!\mathfrak{p}_{H}\rightarrow\mathfrak{k}$ is therefore
injective or zero. Assume (changing notation) that
$\mathfrak{h}_{1}\subset \mathfrak{h}$ is a $\theta$-stable
semi-simple subalgebra such that the injection
$\mathfrak{p}_{1}\subset \mathfrak{p}_{H}$ is holomorphic (for the
choice of $\zeta_{1}\in Z(K_{1})$ where $K_{1}$ is the obvious
maximal compact subgroup of\break
$H_{1}=\exp(\mathfrak{h}_{1})\subset H)$.

It suffices then to check that $F_{c}=0$ on $\mathfrak{p}_{1}$. But any Hermitian Lie
algebra $\mathfrak{h}$ contains a subalgebra $\mathfrak{h}_{1}$ isomorphic to
$\mathfrak{sl}(2,\mathbb{R})$, the injection being holomorphic in the obvious sense
(in fact it contains $\mathfrak{sl}(2,\mathbb{R})^{r}$ where $r$ is the real rank (see
e.g. Ch.~5 of \cite{4}). Thus we are reduced to the case when
$\mathfrak{h}\cong\mathfrak{sl}(2,\mathbb{R})$.

We can also replace $G$ by a larger group. By the results of Satake, $X_{G}$ embeds
into $X_{G_{1}}$ where $G_{1}=SU(p,p)$, as a totally geodesic subvariety, via an
equivariant embedding. Finally we are reduced to the case when $H$ is locally
isomorphic to $SL(2,\mathbb{R})$ or $SU(1,1)$ and $G$ to $SU(p,p)$. (Note that this
does not apply when $G$ has exceptional\break factors).

\section{Computations}

In this paragraph we consider the case, to which we are reduced,
when $H=SU(1,1)$ and $G=SU(p,p)$. We try to solve the linear
algebra problem of \pa 2~--~find $F$ verifying (1)--(3). We have
 \begin{align}
\mathfrak{h}&=\bigg\{\begin{pmatrix} a&z\\\bar{z} &-a
\end{pmatrix}\!:z\in\mathbb{C}\ ,\ a\subset  i\mathbb{R}\bigg\},\\[.4pc]
\mathfrak{g}&=\bigg\{\begin{pmatrix} A&Z\\ {}^{t}{\bar{Z}} &B
\end{pmatrix}\!:\hbox{Tr}(A)+\hbox{Tr}(B)=0\bigg\},
\end{align}
where the block matrices are of size $p\times p$, $Z$ is (complex) arbitrary and $A,B$
are skew-hermitian. Let $u=\left(\begin{smallmatrix}
 &1 \\1 &  \end{smallmatrix}\right)$, $v=\left(\begin{smallmatrix}
 &i \\-i &  \end{smallmatrix}\right)$, $w=\left(\begin{smallmatrix}
 i \\&-i &  \end{smallmatrix}\!\!\right)$, a basis of $\mathfrak{h}$ (the empty entries are zero). Let
 \begin{equation*}
x=\begin{pmatrix}
&1 \\&{}  \end{pmatrix}\quad y=\begin{pmatrix}
&{} \\1  \end{pmatrix}\ ,\ h=\begin{pmatrix}
1\\ &-1  \end{pmatrix}\ ,
\end{equation*}
a basis of $\mathfrak{h}\otimes\mathbb{C}=\mathfrak{sl}(2,\mathbb{C})$. We take
$\mathfrak{k}\subset \mathfrak{g}$ given by block-diagonal matrices, so
\begin{equation*}
\mathfrak{g}=\mathfrak{k}\oplus\mathfrak{p}\ ,\ \mathfrak{p}=\bigg\{\begin{pmatrix}
 &Z \\Z^{\ast}  \end{pmatrix},Z\in M_{p}(\mathbb{C})\bigg\},
\end{equation*}
where $Z^{\ast}={}^{t}\bar{Z}$. Similarly,
\begin{equation*}
\mathfrak{h}=\mathfrak{k}_{H}\oplus\mathfrak{p}_{H}\ ,\ \mathfrak{k}_{H}=\mathbb{R}w\
,\ \mathfrak{p}_{H}=\mathbb{R}v+\mathbb{R}w.
\end{equation*}
If $F\!\!:\mathfrak{h}\rightarrow\mathfrak{g}$ verifies (3) we have
\begin{align}
F(u)&=\begin{pmatrix}
A &Z \\ Z^{\ast}&B  \end{pmatrix},\\[.4pc]
F(v)&=\begin{pmatrix} C &iZ \\-i Z^{\ast}&D  \end{pmatrix},
\end{align}
$A,\dots,D$ verifying of course (5). Let $X,Y,H$ be the images of $x,y,h$. Using (5),
(6) and (7) we have
\begin{align}
X&=\begin{pmatrix}
E &Z \\ {}&F  \end{pmatrix},\\[.4pc]
Y&=\begin{pmatrix}
-E^{\ast}\\ Z^{\ast} &-F^{\ast}  \end{pmatrix},\\[.4pc]
H=[X,Y]&=\begin{pmatrix} -[E,E^{\ast}]+ZZ^{\ast}
&-ZF^{\ast}+E^{\ast}Z\\[.2pc]
FZ^{\ast}-Z^{\ast}E &-[F,F^{\ast}]-Z^{\ast}Z  \end{pmatrix},
\end{align}
where $E,F$ and $Z$ are arbitrary $p\times p$-matrices (with
$\hbox{Tr}(E)+\hbox{Tr}(F)=0$). Since $h=i^{-1}w$, $H$ is block-diagonal by (1);
conjugating $w$ under $K=S(U(p)\times U(p))$ we can assume that the block-diagonal
entries of $H$ are diagonal matrices $H_{1}$, $H_{2}$. The eigenvalues of $H$ are
integral, and constitute the eigenvalues of a representation of
$\mathfrak{sl}(2,\mathbb{C})$.

Let $V\cong\mathbb{C}^{2p}$ be the space of the natural representation of $G$, and
$V=V_{+}\oplus V_{-}$ its decomposition into a positive and a negative subspace.
Then\vspace{-.2pc}
\begin{align*}
&E\!:V_{+}\rightarrow V_{+},\\[.1pc]
&F\!:V_{-}\rightarrow V_{-},\\[.1pc]
&Z\!:V_{-}\rightarrow V_{+}.\\[-1.5pc]
\end{align*}

Let $\lambda_{1}>\cdots>\lambda_{t+1}$ be the distinct eigenvalues of $H$ in $V_{+}$
and $\mu_{1} >\cdots>\mu_{s+1}$ the eigenvalues in $V_{-}\!:s,t\geq 0$. We can write
$V=V^{\mathrm{even}}\oplus V^{\mathrm{odd}}$, the eigenvalues being even or odd in
each summand; this decomposition is preserved by $X$ and $Y$. The decomposition is
orthogonal and compatible with $V=V_{+}\oplus V_{-}$. If $v$ belongs to the
$\lambda$-eigenspace of $V_{+}$ (resp. $V_{-}$), $Ev$ (resp. $Zv, Fv$) belongs to the
$(\lambda+2)$-eigenspace of $V_{+}$ (resp. $V_{+}, V_{-}$)

Consider first the odd part of $V$. We can write in $V_{+}^{\mathrm{odd}}$:
\begin{align*}
E=\begin{pmatrix}
0&E_{1}  \\[.1pc]
&0 &E_{2}\\[.1pc]
&&0 &\ddots\\[.1pc]
&&&\ddots&E_{t}\\[.1pc]
&&&&0
\end{pmatrix}, \quad E^{\ast}=\begin{pmatrix}
0& \\[.1pc]
E_{1}^{\ast}&\ddots\\[.1pc]
 &\ddots &\ddots\\[.1pc]
&&E_{t}^{\ast} &0
\end{pmatrix}.
\end{align*}

Writing $\mathrm{diag}(A_{1},\dots,A_{t+1}) $ for a {\it block}-diagonal matrix we
have
\begin{align*}
EE^{\ast} &=\mathrm{diag}(E_{1}E_{1}^{\ast},\dots,E_{t}E_{t}^{\ast},0)\\[.1pc]
E^{\ast}E &=\mathrm{diag}(0,E_{1}^{\ast}E_{1},\dots,E_{t}^{\ast}E_{t}).
\end{align*}

According to (10),
\begin{align}
-[E,E^{\ast}]+ZZ^{\ast}&=\mathrm{diag}(-E_{1}E_{1}^{\ast},E_{1}^{\ast}E_{1}-E_{2}E_{2}^{\ast},\dots,E_{t}^{\ast}E_{t})+ZZ^{\ast}\nonumber\\[.2pc]
&=\mathrm{diag}(\lambda_{1},\lambda_{2},\dots,\lambda_{t+1}),
\end{align}
where the eigenvalues are now those in $V_{+}^{\mathrm{odd}}$, the last `diagonal'
matrix including of course the multiplicities. Considering the restriction of the
corresponding Hermitian forms to the last summand we see that\vspace{-.3pc}
\begin{equation*}
E_{t}E_{t}^{\ast}+ZZ^{\ast}=\lambda_{t+1}\geqslant 0;
\end{equation*}
since the representation is odd, $\lambda_{1}>\cdots>\lambda_{t+1}>0$.

Similarly in $V_{-}^{\mathrm{odd}}$:\vspace{-.3pc}
\begin{align}
F &=\begin{pmatrix}
0&F_{1} \\[.1pc]  &0 &F_{2}\\[.1pc] &&0\\[.1pc] &&&\ddots &F_{s}\\[.1pc] &&&&0
\end{pmatrix}, \quad F^{\ast}=\begin{pmatrix}
0 \\[.1pc] F_{1}^{\ast} &0\\[.1pc] &&\ddots\\[.1pc] &&F_{s}^{\ast}&0
\end{pmatrix},\nonumber\\[.4pc]
-[F,F^{\ast}]-Z^{\ast}Z
&=\mathrm{diag}(-F_{1}F_{1}^{\ast},\dots,F_{s}^{\ast}F_{s}-Z^{\ast}Z\nonumber\\[.1pc]
&=\mathrm{diag}(\mu_{1},\mu_{2},\dots,\mu_{s+1})
\end{align}
whence $0>\mu_{1}>\cdots>\mu_{s+1}$.\pagebreak

Finally the only non-vanishing part of $Z$ is a map $Z_{1}\!\!:V_{-}(-1)\rightarrow
V_{+}(1)$ (where $V(\lambda)$, $V_{\pm}(\lambda)$ denote the eigenspaces of $H$). Thus
\begin{align}
ZZ^{\ast}&=\mathrm{diag}(0,0,\dots,Z_{1}Z_{1}^{\ast}),\\[.2pc]
Z^{\ast}Z&=\mathrm{diag}(Z_{1}^{\ast}Z_{1},0,\dots,0).
\end{align}

By (11) and (13),
\begin{equation}
\mathrm{diag}(-E_{1}E_{1}^{\ast},E_{1}^{\ast}E_{1}-E_{2}E_{2}^{\ast},\dots,E_{t}^{\ast}E_{t}+Z_{1}Z_{1}^{\ast})
=(\lambda_{1},\dots,\lambda_{t+1})
\end{equation}
with positive eigenvalues. This is impossible unless
\begin{align}
\begin{cases}
t=0 ,\ \lambda_{1}=1 ,\ E=0\ ,\\[.2pc]
Z=Z_{1} ,\ ZZ^{\ast}=1\ .
\end{cases}
\end{align}

(The identity (15) implies that $t=0$; since there is only one
eigenvalue, the representation theory of $SL(2)$ forces it to be
$1$.)

This implies of course that the only eigenvalue $\mu$ is $-1$, so
$s=0$ and $F=0$. Since $E,F$ vanish the embedding is totally
geodesic; the representation of $SL(2)$ or $SU(1,1)$ is a multiple
of the standard representation, in conformity with Satake's
results.

Consider now the even part of $V$. The first part of the argument still applies,
yielding now
\begin{align}
&\lambda_{1}>\cdots>\lambda_{t+1}\geqslant 0,\\[.2pc]
&0\geqslant\mu_{1}>\cdots>\mu_{s+1}.
\end{align}

Now $Z$ is the sum of
\begin{align*}
&Z_{1}\!:V_{-}(0)\rightarrow V_{+}(2),\\[.2pc]
&Z_{2}\!:V_{-}(-2)\rightarrow V_{+}(0).
\end{align*}
Thus
\begin{align}
ZZ^{\ast}&=\mathrm{diag}(0,0,\dots,0,Z_{1}Z_{1}^{\ast},Z_{2}Z_{2}^{\ast}),\\[.2pc]
Z^{\ast}Z&=\mathrm{diag}(Z_{1}^{\ast},Z_{1},Z_{2}^{\ast}Z_{2},0,\dots,0).
\end{align}

By (11) and (19),
\begin{align}
&\hskip -4pc -[E_{1}E^{\ast}]+ZZ^{\ast}\nonumber\\[.3pc]
&\hskip -4pc
\quad=(-E_{1}E_{1}^{\ast},E_{1}^{\ast}E_{1}-E_{2}E_{2}^{\ast},\dots,E_{t-1}^{\ast}E_{t-1}-E_{t}E_{t}^{\ast}
+Z_{1}Z_{1}^{\ast}, E_{t}^{\ast}E_{t}+Z_{2}Z_{2}^{\ast})\nonumber\\[.3pc]
&\hskip -4pc \quad=(\lambda_{1},\dots,\lambda_{t},\lambda_{t+1}),
\end{align}
where we assume so far that both $2$ and $0$ are eigenvalues in
$V_{+}^{\mathrm{even}}$. This implies first that there are only two eigenvalues since
$-E_{1}E_{1}^{\ast}=\lambda_{1}>0$ for $t>1$. Furthermore, the last entry in (21)
yields $E_{1}E_{1}^{\ast}+Z_{2}Z_{2}^{\ast}=0$, whence $E=E_{1}=0$ and $Z_{2}=0$.

If $2$ does not occur in $V_{+}^{\mathrm{even}}$, the representation on
$V^{\mathrm{even}}$ is trivial; if $0$ does not occur $Z_{2}$ is absent. In this case,
\begin{align}
ZZ^{\ast} &=\mathrm{diag}(0,\dots 0,Z_{1}Z_{1}^{\ast}),\nonumber\\[.2pc]
Z^{\ast}Z &=\mathrm{diag}(Z_{1}^{\ast}Z_{1},0,\dots 0)
\end{align}
and
\begin{align*}
-[E,E^{\ast}]+ZZ^{\ast}&=(-E_{1}E_{1}^{\ast},E_{1}^{\ast}E_{1}-E_{2}E_{2}^{\ast},\dots,
E_{t}^{\ast}E_{t}+Z_{1}Z_{1}^{\ast})\\[.2pc]
&=(\lambda_{1},\dots,\lambda_{t+1})
\end{align*}
with $\lambda_{t+1}=2$. This equality $(-E_{1}E_{1}^{\ast}=\lambda_{1}>0)$ implies
that there is only one eigenvalue $(t=0)$ and therefore $E=0$.

Of course, a similar computation, as in the odd case, applies to the negative part,
using now (20); if there are two eigenvalues $(0,-2)$ we deduce that
\begin{equation*}
F=F_{1}=0\ ,\ Z_{1}=0.
\end{equation*}

Thus $Z=0$, contrary to the assumption that it represented the tangent map to an equivariant embedding.

Finally, consider the case where $0$ does not occur in $V_{+}^{\mathrm{even}}$ or
$V_{-}^{\mathrm{even}}$.

The computations being symmetric we can assume for instance that it is missing in
$V_{+}^{\mathrm{even}}$; we already know that the eigenvalue $2$ only occurs, so the
eigenvalues in $V^{\mathrm{even}}$ are $(2,0,-2)$; moreover $E=F=0$ by the arguments
given already, so the embedding should be totally geodesic. We know that this is
impossible, by Satake's results. In fact, using (22) and (12) we see that
\begin{equation*}
\mathrm{diag}(-F_{1}F_{1}^{\ast}-Z_{1}Z_{1}^{\ast},F_{1}^{\ast}F_{1})=(\mu_{1},\mu_{2})=(0,-2)
\end{equation*}
which is impossible.


\begin{thebibliography}{9}
\bibitem{1} Andr\'e~Y, G-functions and geometry, {\it Aspects of Math.}
Vieweg (ed.) (1989)

\bibitem{2} Chen~B~Y and Nagano~T, Totally geodesic
submanifolds of symmetric spaces~I, {\it Duke Math. J.} {\bf 46} (1977) 745--755

\bibitem{3} Ihara~S-I, Holomorphic imbeddings of symmetric domains,
{\it J. Math. Soc. Japan} {\bf 19} (1967) 261--302; Suppl.
543--544

\bibitem{4} Mok~N, Metric rigidity theorems on hermitian locally
symmetric manifolds (Singapore: World Scientific) (1989)

\bibitem{5} Moonen~B, Linearity properties of Shimura varieties~I,
{\it J. Alg. Geom.} {\bf 7} (1998) 539--567

\bibitem{6} Moret-Bailly~L, Groupes de Picard et
probl\`emes de Skolem~II, {\it Ann. Sc. Ecole Normale Sup.~(4)} {\bf 22} (1989)
181--194

\bibitem{7} Oort~F, Canonical liftings and dense sets of CM-points,
in: Arithmetic geometry (Cortona) (1994), Symp. Math. XXXVII
(Cambridge: Cambridge Univ. Press) (1997)

\bibitem{8} Satake~J, Holomorphic imbeddings of symmetric domains
into a Siegel space, {\it Am. J. Math.} {\bf 87} (1965) 425--461
\end{thebibliography}
\end{document}